\documentclass[a4paper]{article}
\usepackage{delarray,verbatim,enumerate}
\usepackage{lscape}
\usepackage{pdflscape}
\usepackage{amsmath,amsthm,amstext,amsbsy,amssymb, amsfonts,amscd}

\newcommand{\address}[1]{\vskip3mm\noindent#1}

\newtheorem{proposition}{Proposition}
\newtheorem{definition}{Definition}
\newtheorem{theorem}{Theorem}

\newtheorem{remark}{Remark}
\def\goth{\mathfrak}
\newcommand{\Z}        {{\mathbb Z}}

\newcommand{\Fq}       {{\mathbb F}_{2}}
\newcommand{\Q}        {{\mathbb Q}}

\newcommand{\RR}       {{\cal R}}
\def\wi{\widetilde}		\def\rg{\operatorname{rk}}	
\def\div{\operatorname{div}}
\def\sg{\operatorname{sg}}
\def\Gal{\operatorname{Gal}}	
	
\newcommand {\IDa} {{\goth a}}

\newcommand {\IDb} {{\goth b}}
\newcommand {\IDc} {{\bf c}}
\newcommand {\IDe} {{\bf e}}
\newcommand {\IDf} {{\bf f}}
\newcommand {\IDg} {{\bf g}}

\newcommand {\IDp} {{\goth p}}

\newcommand {\IDq} {{\goth q}}

\newcommand{\Cl}		{{\mathcal C \ell}}
\newcommand{\Dl}		{{\mathcal D \ell}}
\newcommand{\Pl}		{{\mathcal P \ell}}

\title{\LARGE \bf  Computation of 2-groups of positive classes 
 of exceptional number fields\footnote{Version de travail du 03 decembre 2007}}
\author{Jean-Fran\c cois {\sc Jaulent}, Sebastian {\sc Pauli}, \\
 Michael E. {\sc Pohst} \& Florence {\sc Soriano--Gafiuk}}

\date{}
\def\LaTeX{L\kern-.25em\raise.425ex\hbox{a}\kern-.075em\TeX}

\begin{document}
\label{firstpage}
\maketitle
%

\bigskip\bigskip

{\small
\noindent {\bf Abstract.} We present an algorithm for computing 
the 2-group $\Cl_F^{\,pos}$ of the positive divisor classes in case 
the number field $F$ has exceptional dyadic places. As an application, 
we compute the 2-rank of the wild kernel $W\!K_2(F)$ in $K_2(F)$. 
\medskip

\noindent {\bf R\'esum\'e}. Nous d\'eveloppons un algorithme pour 
d\'eterminer le 2-groupe $\Cl_F^{\,pos}$ des classes positives dans 
le cas o\`u le corps de nombres consid\'er\'e $F$ poss\`ede des places 
paires exceptionelles. Cela donne en particulier
le 2-rang du noyau sauvage $W\!K_2(F)$.
}


\section{Introduction}


The logarithmic $\ell$-class group $\wi\Cl_F$ whas introduced in \cite{J2} by J.-F. Jaulent who used it to study the $\ell$-part $W\!K_2(F)$ of the wild kernel in number fields: if $F$ contains a primitive $2\ell^t$-th root of unity ($t>0$), there is a natural isomorphism\smallskip

\centerline{$\mu_{\ell^t} \otimes_\Z \wi\Cl_F \simeq W\!K_2(F)/W\!K_2(F)^{\ell^t}$,}\smallskip

\noindent so the $\ell$-rank of $W\!K_2(F)$ coincides with the $\ell$-rank of the logarithmic  group $\wi{\Cl}_F$.
An algorithm for computing $\wi{\Cl}_F$ for Galois extensions $F$ was developed in \cite{DS} and later generalized and improved for arbitrary number fields in \cite{DJ+}.\smallskip

In case the prime $\ell$ is odd, the assumption $\mu_\ell \subset F$ may be easily passed if one considers the cyclotomic extension $F(\mu_{\ell})$ and gets back to $F$ via the so-called transfer (see \cite{JS1}, \cite{PS} and \cite{So2}). However for $\ell=2$ the connection between symbols and logarithmic classes is more intricate: in the non-exceptional situation ({\it i.e.} when the cyclotomic $\Z_2$-extension  $F^c$ contains the fourth root of unity $i$) the 2-rank of $W\!K_2(F)$ still coincides with the 2-rank of  $\wi\Cl_F$. Even more if the number field $F$ has no exceptional dyadic place ({\it i.e.} if one has $i \in F^c_\IDq$ for any $\IDq | 2$),  the same result holds if one replace the ordinary logarithmic class group $\wi\Cl_F$ by a narrow version $\wi\Cl{}_F^{\, res}$. The algorithmic aspect of this is treated in \cite{JPPS}.\smallskip

Last in \cite{JS2} the authors pass the difficulty in the remaining case by introducing  a new 2-class groups $\Cl_F^{\,pos}$, the {\it 2-group of positive divisor classes}, which satisfies the rank identity:
$\rg _2 \Cl^{\,pos}_{F} = \rg_2 W\!K_2(F)$.\smallskip

In this paper we develop an algorithm for computing both $\Cl_F^{\,pos}$  and  $\wi\Cl{}_F^{\,pos}$ in case the number field $F$ does contain exceptional dyadic places.\nopagebreak

We conclude with several examples.  Combining our algorithm with the work of Belabas and Gangl \cite{BG} on the computation of  the tame kernel of $K_2$ we obtain the complete structure of the wild kernel in some cases.


\section{Positive divisor classes of degree zero}

 
\subsection{The group of logarithmic divisor classes of degree zero}

 
Throughout this paper the prime
number $\ell$ equals 2 and we let $i$ be a primitive fourth root
of unity. Let $F$ be a number field of degree $n=r+2c$.  
According to \cite{J1}, for every place $\IDp$ of $F$
there exists a 2-adic valuation $\widetilde{v}_{\IDp}$ which is
related to the wild 2-symbol in case the cyclotomic 
$\Z_2$-extension of $F_{\IDp}$ contains $i$. The degree $\deg \IDp$ of
$\IDp$ is a 2-adic integer such that the image of the map Log$\;
|\;|_{\IDp}$ is the $\Z_2$-module $\deg (\IDp) \; \Z_2$ (see
\cite{J2}). (By Log we mean the usual $2$-adic logarithm.) The
construction of the 2-adic logarithmic valuations
$\widetilde{v}_{\IDp}$ yields
\begin{equation}
 \forall \alpha \in  \RR_F:=\Z_2 \otimes_{\Z} F^{\times} \; :
 \; \sum_{\IDp \in Pl^{\,0}_F} \widetilde{v}_{\IDp}(\alpha) \deg (\IDp) \; =\; 0,
\end{equation}
where $Pl^{\,0}_F$ denotes the set of finite places of the number field $F$. Setting
\[
\wi\div(\alpha) \; :=\; \sum_{\IDp \in Pl^{\,0}_F} \wi{v}_{\IDp}(\alpha)
 \IDp
\]
we obtain by $\Z_2$-linearity:
\begin{equation}
\deg (\wi\div(\alpha)) \; =\; 0 .
\label{degdiv}
\end{equation}
We define the 2-{\em group of logarithmic divisors of degree 0} as the kernel of the degree map 
$\deg$  in the direct sum $\Dl_F= \sum_{\IDp \in Pl^{\,0}_F} \Z_2 \; \IDp$:
\[
\widetilde{\Dl}_F \; :=\; \;\left\{ \textstyle\sum_{\IDp \in Pl^{\,0}_F} a_{\IDp} \IDp \in
\Dl _F \mid  {\textstyle \sum_{\IDp \in Pl^{\,0}_F} a_{\IDp}} \deg (\IDp)
=0 \right\};
\]
and the {\em subgroup of principal logarithmic divisors} as the image of the logarithmical map $\wi\div$:
\[
\wi{\Pl}_F \; :=\; \{ \widetilde\div(\alpha) \mid \alpha \in \RR_F
 \} \;\; .
\]
Because of (\ref{degdiv}) $\wi{\Pl}_F$ is clearly a subgroup of $\wi{\Dl}_F$.
More ever by the so-called extended Gross conjecture, the factorgroup
\[
\wi{\Cl}_F \; :=\; \widetilde{\Dl}_F / \widetilde{\Pl}_F
\]
is a finite 2-group, the {\em 2-group of logarithmic divisor classes}. So, under this conjecture,
$\wi\Cl_F$ is just the torsion subgroup of the group $$\Cl_F:=\Dl_F/\wi\Pl_F$$ of logarithmic classes
(without any asumption of degree). 

\begin{remark}
{\rm Let $F^{+}$ be the set of all totally positive elements of $F^\times$ ({\it i.e.} the subgroup $F^+:=\{x \in F^\times |\;x_\IDp > 0$ \quad  for all real $\IDp\}$.  For
\[
\widetilde{\Pl}{}_F^{+} \; :=\; \{ \widetilde\div(\alpha) \mid \alpha \in
\RR_F^{+}:=\Z_2 \otimes_\Z F^{+} \}
\]
the factor group
\[
\Cl_F^{\,res} \; :=\; \Dl_F / \wi\Pl{}_F^{+} \qquad {\rm (resp.}\quad \wi\Cl{}_F^{\,res} \; :=\; \wi{\Dl}_F / \wi{\Pl}{}_F^{+})
\]
is the {\em 2-group of narrow logarithmic divisor classes} of the number field $F$ (resp. the {\em 2-group of narrow logarithmic divisor classes of degree 0}) introduced in \cite{So1} and computed in \cite{JPPS}.}
\end{remark}


\subsection{Signs and places}


For a field $F$ we denote by $F^c$, (respectively $F^c[i]$) the cyclotomic $\Z_2$-extension (resp. the maximal cyclotomic pro-2-extension) of $F$.

\vspace{0.3cm}
We adopt the notations and definitions in this section from \cite{JS2}.

\begin{definition}[{\em signed places}] \rm
Let $F$ be a number field.
We say that a non-complex place $\IDp$ of $F$ is {\em signed} if and only if $F_{\IDp}$ does not contains the fourth root $i$. These are the places which do not decompose in the extension $F[i]/F$. \\
We say that $\IDp$ is {\em logarithmically signed} if and only if the cyclotomic $\Z_2$-extension $F_{\IDp}^c$ does not contain $i$. These are the places which do not decompose in $F^c[i]/F^c$.
\end{definition}

\begin{definition}[{\bf sets of signed places}] \rm
By $P\!S$, respectively $P\!L\!S$, we denote the sets of signed, respectively logarithmically signed, places:
\begin{eqnarray*}
P\!S & := & \{ \IDp \mid i \not\in F_{\IDp} \}  \;\; , \\
P\!L\!S & := & \{ \IDp \mid i \not\in F_{\IDp}^c \}  \;\; .
\end{eqnarray*}
A finite place $\IDp \in P\!L\!S$ is called {\em exceptional}. The set of exceptional places is denoted by $P\!E$. Exceptional places are even (i.e. finite places dividing 2).
\end{definition}

These sets satisfy the following inclusions:
\[
P\! \subset P\!L\!S  = P\!E \cup P\!R \subset Pl(2) \cup Pl(\infty)
\]
where $Pl(2),\, Pl(\infty), \, P\!R$ denote the sets of even, infinite and real places of $F$, respectively. From this the finiteness of $P\!L\!S$ is obvious.\smallskip

We recall the canonical decomposition $\Q_2^{\times}=2^{\Z} \times (1+4\Z_2) \times \langle -1
\rangle$ and we denote by $\epsilon$ the projection from $\Q_2^{\times}$
onto $\langle -1 \rangle$.

\begin{definition}[{\bf sign function}] \rm
For all places $\IDp$ we define a sign function via
\[
\sg_{\IDp} \; :\; F_{\IDp}^{\times} \rightarrow \langle -1 \rangle
\; :\; x \mapsto 
\left\{ \begin{array}{lll}
1 & \mbox{for} & \IDp \; \mbox{complex} \\
\mbox{sign}(x) & \mbox{for} & \IDp \; \mbox{real} \\
\epsilon (N \IDp^{-\nu_{\IDp}(x)}) & \mbox{for} & \IDp \not\,\mid 2\infty \\
\epsilon (N_{K_{\IDp} / \Q_{2}} (x) N \IDp^{-\nu_{\IDp}(x)} ) &
\mbox{for} & \IDp \mid 2 
\end{array} \right. \;\; .
\]
These sign functions satisfy the product formula:
\[
\forall x \in F^\times \qquad \prod_{\IDp \in Pl_F} \sg (x) =1.
\]
\end{definition}

In addition we have:

\begin{proposition}
The places $\IDp$ of $F$ satisfy the following properties:
\begin{enumerate}
\item[(i)] if $\IDp \in P\!L\!S$ then $(\sg_{\IDp}, \widetilde{v}_{\IDp})$ is
surjective; 
\item[(ii)] if $\IDp \in P\!S \setminus P\!L\!S$ then
$\sg_{\IDp}(\;\;) =(-1)^{\widetilde{v}_{\IDp}(\;\;)}$ and
$\widetilde{v}_{\IDp}$ is surjective; 
\item[(iii)] if $\IDp \not\in P\!S$ then
$\sg_{\IDp} (F_{\IDp}^{\times})=1$ and $\widetilde{v}_{\IDp}$ is surjective.
\end{enumerate}
\end{proposition}

\begin{remark}\rm
The logarithmic valuation $\widetilde{v}_{\IDp}$ is surjective in all three cases.
Part 2 of the preceding result is often used for testing $\IDp \in P\!L\!S$.
\end{remark}


\subsection{The group of positive divisor classes}


For the introduction of that group we modify several
notations from \cite{JS2} in order to make them suitable
for actual computations.\smallskip

Since $P\!L\!S$ is finite we can fix the order of the logarithmically signed places,
say $P\!L\!S=\{\IDp_1,\cdots,\IDp_m\}$, with $P\!E=\{\IDp_1,\cdots,\IDp_e\}$  and $P\!R=\{\IDp_{e+1},\cdots,\IDp_m\}$.
Accordingly we define vectors $\IDe=(e_1,\cdots,e_m) \in \{ \pm 1\}^m$.\smallskip

For each divisor $\IDa = \sum_{\IDp \in Pl^{\,0}_F} a_{\IDp} \IDp$, we form pairs $(\IDa,\IDe)$ and put
\begin{equation}
 \sg (\IDa,\IDe) \; :=\; \prod_{\IDp \in P\!S \setminus P\!L\!S} 
\, (-1)^{a_{\IDp}}  \times \prod_{i=1}^m e_i 
\label{defd}
\end{equation}
Let $\Dl_F (P\!E):=\left\{\IDa\in \Dl_F\,\big|\,\IDa=\sum_{\IDp \in P\!E} a_{\IDp} \IDp\right\}$ be the $\Z_2$-submodule of $\Dl_F$ generated by the exceptional dyadic places. And let $\Dl^{P\!E}_F$ be the factor group $\Dl_F/\Dl_F(P\!E)$. Thus the {\em group of positive divisors} is the $\Z_2$-module:
\begin{equation}
\Dl_F^{\,pos} \; :=\; \left\{ (\IDa,\IDe) \in \Dl_F^{P\!E} \times \{ \pm 1 \}^m \,\Big|\,  \sg (\IDa,\IDe)=1 \right\}
\label{defdlpos}
\end{equation}
For $\alpha \in \RR_F:=\Z_2 \otimes_\Z F^\times$, let $\wi\div{}'(\alpha)$ denotes the image of $\wi\div(\alpha)$ in $\Dl_F^{P\!E}$ and $\sg (\alpha )$ the vector of  signs $(\sg_{\IDp_1}(\alpha),\dots,\sg_{\IDp_m}(\alpha))$ in $\{\pm1\}^m$.
Then 
\begin{equation}
\widetilde{\Pl}{}_F^{\,pos} \; :=\; \left\{ (\wi\div{}'(\alpha), 
\sg (\alpha ))  \in \Dl_F^{P\!E} \times \{ \pm 1\}^m \,\Big|\,
\alpha \in \RR_F \right\} 
\label{defplpostil}
\end{equation}
is obviously a submodule of $\Dl_F^{\,pos}$ which is called the {\em principal submodule}.

\begin{definition}[{\bf positive  divisor classes}] \rm
With the notations above:
\begin{itemize}
\item[(i)] The group of {\em positive logarithmic divisor classes} is the factor group
\[
\Cl_F^{\,pos} \; =\; \Dl_F^{\,pos} / \widetilde{\Pl}{}_F^{\,pos} \;\; .
\]
\item[(ii)] The subgroup of {\em positive logarithmic divisor classes of degree zero} is the kernel $\wi{\Cl}{}_F^{\,pos}$ of the degree map $\deg$ in $\Cl_F^{\,pos}$:
\[
\wi\Cl{}_F^{\,pos}:= \{ (\IDa,\IDe)+\wi\Pl_F^{\,pos} \; | \;  \deg (\IDa) \in \deg (\Dl_F(P\!E)) \} .
\]
\end{itemize}
\end{definition}

\begin{remark}\rm
 The group $\Cl_F^{\,pos}$ is infinite whenever the number field $F$ has no exceptional places, since in this case $\deg(\Cl_F^{\,pos})$ is isomorphic to $\Z_2$. The finiteness of $\Cl_F^{\,pos}$ in case $P\!E \ne \emptyset$ follows from the so-called generalized Gross conjecture.
\end{remark}

For the computation of $\wi\Cl{}_F^{\,pos}$ we need to introduce {\em primitive divisors}.

\begin{definition} \rm
A divisor $\IDb$ of $F$ is called a {\em primitive} divisor if $\deg (\IDb)$ generates
the $\Z_2$-module $\deg (\Dl_F) =4[F \cap \Q^c : \Q]\Z_2$.
\end{definition}

We close this section by presenting a method for exhibiting such a divisor:\smallskip

Let $\IDq_{1},\cdots, \IDq_{s}$ be all dyadic primes; and  $\IDp_{1},\cdots, \IDp_{s}$ be a finite set of non-dyadic primes which generates the 2-group of 2-ideal-classes ${\Cl'_F}$ ({\it i.e.} the quotient of the  usual 2-class group by the subgroup generated by ideals above 2). 

Then every $\IDp \in \{\IDq_1,\cdots,\IDq_s,\IDp_1,\cdots,\IDp_t\}$ with minimal 2-valuation $\nu_2(\deg\IDp)$ is primitive.


\subsection {Galois interpretations and applications to $K$-theory}


Let $F^{lc}$ be  the locally cyclototomic 2-extension of $F$ ({\it i.e.} the maximal abelian pro-2-extension of $F$ which is completely split at every place over the cyclotomic $\Z_2$-extension $F^c$. Then by $\ell$-adic class field theory ({\it cf.} \cite{J1}), one has the following interpretations of the logarithmic class groups:
\[
\Gal (F^{lc}/F) \simeq \Cl_F \qquad {\rm and} \qquad \Gal (F^{lc}/F^c) \simeq \wi\Cl_F\, .
\]

\begin{remark}\rm
Let us assume $i \notin F^c$. Thus we may list the following special cases:
\begin{enumerate}
\item[$(i)$] In case $P\!L\!S=\emptyset$, the group $\Cl^{\,pos}_F \simeq \Z_2 \oplus \wi\Cl{}^{\, pos}_F$ of positive divisor classes has index 2 in the group $\Cl_F \simeq \Z_2 \oplus \wi\Cl_F$ of logarithmic classes 
of arbitrary degree; as a consequence its torsion subgroup $\;\wi\Cl{}^{\, pos}_F$ has index 2 in the finite 
group $\wi\Cl_F$ of logarithmic classes of degree 0 yet computed in \cite{DJ+}.

\item[$(ii)$] In case $P\!E=\emptyset$, the group $\Cl^{\,pos}_F \simeq \Z_2 \oplus \wi\Cl{}^{\, pos}_F$ has index 2 in the group $\Cl^{\,res}_F \simeq \Z_2 \oplus \wi\Cl{}^{\, res}_F$ of narrow logarithmic classes 
of arbitrary degree; and its torsion subgroup $\wi\Cl{}^{\, pos}_F$ has index 2 in the finite 
group $\wi\Cl{}^{\, res}_F$ of narrow logarithmic classes of degree 0 introduced in \cite{So1} and computed in \cite{JPPS}.
\end{enumerate}
\end{remark}

\begin{definition}\rm We adopt the following conventions from \cite{H1,H2,JS2,JS3}:
\begin{enumerate}
\item[$(i)$] $F$ is {\em exceptional} whenever one has $i \notin F^c$ ({\it i.e.} $[F^c[i]:F^c]=2$);
\item[$(ii)$] $F$ is {\em logarithmically signed} whenever one has $i \notin F^{lc}$ ({\it i.e.} $P\!L\!S \ne \emptyset$);
\item[$(iii)$] $F$ is {\em primitive} whenever one at least between the exceptional places does not split in (the first step of the cyclotomic $\Z_2$-extension) $F^c/F$.
\end{enumerate}
\end{definition}\smallskip

The following theorem is a consequence of the results in \cite{H1, H2, J1, J2, JS2, JS3}:\smallskip

\begin{theorem}
Let $W\!K_2(F)$ (resp. $K_2^\infty(F): =\cap_{n \ge 1}K_2^{2^n}(F)$) be be the 2-part of the wild kernel (resp. the 2-subgroup of infinite height elements) in $K_2(F)$.
\begin{enumerate}
\item[$(i)$] In case $i  \in F^{lc}$ (i.e. in case $P\!L\!S = \emptyset$), we have both:
\[
\rg_2 W\!K_2 (F) \, =\, \rg_2 \wi\Cl_F \, = \,  \rg_2 \wi\Cl{}_F^{\,res} .
\]

\item[$(ii)$] In case $i \notin F^{lc}$ but $F$ has no exceptional places (i.e. $P\!E = \emptyset$), we have:
\[
\rg_2 W\!K_2 (F) \, = \,  \rg_2 \wi\Cl{}_F^{\,res} .
\]

\item[$(iii)$]In case $P\!E \neq \emptyset$, then we have
\[
\rg_2 W\!K_2 (F) \, =\, \rg_2 \Cl_F^{\,pos}  .
\]
And in this last situation there are two subcases:
\begin{enumerate}
\item If $F$ is primitive, i.e. if the set $P\!E$ of exceptional  dyadic
places contains a primitive place, we have:
\[
K_2^{\infty} (F) \, =\, W\!K_2(F) \;\; .
\]
\item If $F$ is imprimitive and $K_2^{\infty} (F) =\oplus_{i=1}^n \, \Z /2^{n_i} \Z $, we get:
\begin{enumerate}
\item $W\!K_2(F)  = \Z/2^{n_1+1} \Z \oplus \left( 
\oplus_{i=2}^n \Z/2^{n_i} \Z \right) 
\mbox{ if }  \rg_2 (\wi\Cl{}_F^{\,pos})  = \rg_2 (\Cl_F^{\,pos})$;
\item  $W\!K_2(F)  = \Z/2\Z \oplus \left(
\oplus_{i=1}^n \Z/2^{n_i} \Z \right) 
\mbox{ if } \rg_2 (\wi\Cl{}_F^{\,pos})  < \rg_2 (\Cl_F^{\,pos})$.
\end{enumerate}
\end{enumerate}
\end{enumerate}
\end{theorem}

\newpage

\section {Computation of positive divisor classes}


We assume in the following that the set $P\!E$ of exceptional places is not empty.


\subsection{Computation of exceptional units}


Classically the group of logarithmic units is the kernel in ${\cal R}_F$ of the logarithmic valuations (see \cite{J1}):
\[
\wi{\cal E}_F = \{ x \in \RR_F \mid \forall \IDp : \widetilde{v}_{\IDp}(x)=0 \}
\]
In order to compute positive divisor classes in case  $P\!E$  is not empty, we ought to introduce a new group of units:

\begin{definition} \rm
We define the group of {\em logarithmic exceptional units} as the kernel of the non-exceptional logarihtmic valuations:
\begin{equation}
\wi{\cal E}{}^{e\!x\!c}_F = \{ x \in \RR_F \mid \forall \IDp \notin P\!E : \widetilde{v}_{\IDp}(x)=0 \}
\label{tildeU}
\end{equation}
\end{definition}

We just know that  is a subgroup of the 2-group of 2-units ${\cal E}'_F =\ Z_2 \otimes E'_F$. If we assume that there are exactly $s$ places in $F$ containing $2$ we have, say:
\[
E'_F \; =\; \mu_F \times \langle \varepsilon_1,\cdots, \varepsilon_{r+c-1+s} \rangle
\]
For the calculation of $\wi{\cal E}{}^{e\!x\!c}_F$ we use the same precision $\eta$ as for our $2$-adic approximations used in the course of the calculation of $\,\wi{\Cl}_F$.
Then we obtain a system of generators of $\wi{\cal E}{}^{e\!x\!c}_F$ by computing the nullspace of the matrix
\[
 B \; =\; \left( \begin{array}{cccccc}
   &   & | & 2^\eta &  \cdots & 0 \\
   & \widetilde{v}_{\IDp_i}(\varepsilon_j) & | & \cdot & \cdots & \cdot \\
   &   & | & 0 & \cdots & 2^\eta
\end{array} \right) 
\]
with $r+c-1+s+e$ columns and $e$ rows, where $e$ is the cardinality of $P\!E$ and the precision $\eta$ is determined as explained in \cite{DJ+}. \smallskip

We assume that the nullspace is generated by the columns of the matrix

\[
 B' \; =\; \left( \begin{array}{ccc}
   &   &  \\
   & C &  \\
   &   &  \\
 - & - & - \\
   &   &   \\
   & D &  \\
   &   & 
\end{array} \right) 
\]\bigskip

\noindent where $C$ has $r+c-1+s$ and $D$ exactly $e$ rows. It suffices to consider $C$. Each column $(n_1,\cdots,n_{r+c-1+s})^{tr}$ of $C$ corresponds to a unit
\[
\prod_{i=1}^{r+c-1+s} \varepsilon_i^{n_i} \in \wi{\cal E}{}^{e\!x\!c}_F
\RR_F^{2^\eta}
\]
so that we can choose
\[
\widetilde{\varepsilon}:=\prod_{i=1}^{r+c-1+s} \varepsilon_i^{n_i} 
\]
as an approximation for an exceptional unit. This procedure yields $k \ge r+c+e$ exceptional units, say: $\wi\varepsilon_1,\cdots,\wi\varepsilon_k$.
By the so-called generalized conjecture of Gross we would have exactly $r+c+e$ such units. So we assume in the following that the procedure does give $k=r+c+e$ (otherwise we would refute the conjecture).
Hence, from now on we may assume that we have determined exactly
$r+c+e$ generators $\wi{\varepsilon}_1,\cdots,\widetilde{\varepsilon}_{r+c+1}$ of $\wi{\cal E}{}^{e\!x\!c}_F$, and we write:
\[
\wi{\cal E}{}^{e\!x\!c}_F \; =\; \langle -1 \rangle \times \langle \wi\varepsilon_1,\cdots, \wi\varepsilon_{r+c-1+e} \rangle
\]

\begin{definition} \rm
The kernel of the canonical map ${\cal R}_F \rightarrow \Dl_F^{\,pos}$ is the subgroup of {\em positive logarithmic units}:
\[
\wi{\cal E}{}^{\,pos}_F \; =\; \{ \wi\varepsilon \in \wi{\cal E}{}^{e\!x\!c}_F \; |\; \forall \IDp \in P\!L\!S \quad \sg_\IDp (\wi\varepsilon) =+1 \}
\]

The subgroup  $\wi{\cal E}{}^{\,pos}_F$ has finite index in the group $\wi{\cal E}{}^{e\!x\!c}_F$ of exceptional units.
\end{definition}


\subsection{The algorithm for computing $\Cl_F^{\,pos}$}


We assume  $P\!E \neq \emptyset$ and that the logarithmic 
2-class group $\widetilde{\Cl}_F$ is isomorphic to the direct sum
\[
\widetilde{\Cl}_F \; \cong\; \oplus_{i=1}^{\nu} \, \Z / 2^{n_i}\Z
\]
subject to $1 \leq n_1 \leq \cdots\leq n_{\nu}$. Let $\IDa_i \;\;
(1 \leq i \leq \nu)$
be fixed representatives of the $\nu$ generating
divisor classes. Then any divisor
$\IDa$ of $\Dl_F$ can be written as
\[
\IDa \; =\; \sum_{i=1}^{\nu} \, a_i \IDa_i + \lambda \IDb
+ \widetilde{\div}(\alpha) 
\]
with suitable integers $a_i \in \Z_2$, a primitive divisor
$\IDb$, $\lambda=\frac{\deg (\IDa)}{\deg (\IDb)}$ and an
appropriate element $\alpha$ of $\RR_F$. With
each divisor $\IDa_i$ we associate a vector
\[
\IDe_i \; :=\; (\sg (\IDa_i,{\bf 1}),1,\cdots,1) \in \{ \pm 1 \}^m \;\; ,
\]
where $m$ again denotes the number of divisors in $P\!L\!S$. 
Clearly, that representation then satisfies 
$ \sg (\IDa_i,\IDe_i)=1$, hence the element
$(\IDa_i,\IDe_i)$ belongs to $\Dl_F^{\,pos}$. Setting
$\IDe_{\IDb}= (\sg (\IDb,{\bf 1}),1,\cdots,1)$ as above and
writing
\[
\IDe'\; :=\; \sg (\alpha ) \times
\prod_{i=1}^{\nu} \IDe_i^{a_i} \times \IDe \times \IDe_{\IDb}^{\lambda}
\]
for abbreviation any element $(\IDa, \IDe)$ of $\Dl_F^{\,pos}$
can then be written in the form
\begin{eqnarray*}
(\IDa,\IDe) & = & \left( \sum_{i=1}^{\nu} \, a_i \IDa_i + \lambda
\IDb + \widetilde{\div}(\alpha),  \IDe'
 \times \prod_{i=1}^{\nu} \IDe_i^{a_i} \times \sg (\alpha)
\times  \IDe_{\IDb}^{\lambda} \right) \\
 & = & \sum_{i=1}^{\nu} \, a_i (\IDa_i,\IDe_i) + 
\lambda (\IDb,\IDe_{\IDb}) +
({\bf 0}, \IDe') + (\widetilde{\div}(\alpha),\sg (\alpha )) \;\; .
\end{eqnarray*}
The multiplications are carried out coordinatewise. 
The vector $\IDe'$ is 
therefore contained in the $\Z_2$-module generated by
$\IDg_i \in \Z^m \;\; (1 \leq i \leq m)$ with 
$\IDg_1=(1,\cdots,1)$,
whereas $\IDg_i$ has first and $i$-th coordinate -1,
all other coordinates 1 for $i >1$.

As a consequence, the set
\[
\{ (\IDa_j,\IDe_j) \mid 1 \leq j \leq \nu \} \cup \{ (0,\IDg_i) \mid
2\leq i \leq m \} \cup \{ (\IDb,\IDe \}
\]
contains a system of generators of
$\Cl_F^{\,pos}$  ( note that $(0,\IDg_1)$ is trivial in
$\Cl_F^{\,pos}$).

\noindent We still need
to expose the relations among those. But the latter are easy to
characterize. We must have
\begin{eqnarray*}
 \sum_{j=1}^{\nu} \, a_j (\IDa_j,\IDe_j) + \sum_{i=2}^m
b_i ({\bf 0},\IDg_i) + \lambda (\IDb,\IDe_{\IDb})
 & \equiv & 0 \, \bmod \; \wi{\Pl}{}_F^{\,pos} \;\; ,\\
 \sum_{j=1}^{\nu} \, a_j (\IDa_j,\IDe_j) + \sum_{i=2}^m
b_i ({\bf 0},\IDg_i)  + \lambda (\IDb,\IDe_{\IDb}) & = &
(\widetilde{\div}(\alpha), \sg (\alpha )) + \sum_{\IDp \in P\!E}
 \; (d_{\IDp} \IDp, {\bf 1})
\end{eqnarray*}
with indeterminates $a_j,b_i,d_{\IDp}$ from $\Z_2$. Considering the
two components separately, we obtain the conditions
\begin{equation}
\sum_{j=1}^{\nu} \, a_j \IDa_j  + \lambda \IDb\; \equiv \;
\sum_{\IDp \in P\!E} \; d_{\IDp} \IDp \, \bmod \; \widetilde{\Pl}_F
\label{cond1}
\end{equation}
and
\begin{equation}
\prod_{j=1}^{\nu} \IDe_j ^{a_j} \times \prod_{i=2}^m \IDg_i^{b_i} 
 \times \IDe_{\IDb}^{\lambda} \; =\;  \sg (\alpha ) \;\; .
\label{cond2}
\end{equation}
Let us recall that we have already ordered $P\!L\!S$ so that exactly the first $e$ elements
$\IDp_1,\cdots,\IDp_{e}$ belong to $P\!E$.
Then the first one of the conditions above
is tantamount to
\[
\sum_{j=1}^{\nu} \, a_j \IDa_j \; \equiv \;
\sum_{i=1}^{e} \, d_{\IDp_i} \left( \IDp_i -
\frac{\deg \IDp_i}{\deg \IDb} \IDb \right) \; \bmod \; \widetilde{\Pl}_F
\;\; .
\]
The divisors
\[
 \IDp_i - \frac{\deg \IDp_i}{\deg \IDb} \IDb
\]
on the right-hand side can again be expressed by the
$\IDa_j$. For $1 \leq i \leq e$ we let
\[
\widetilde{\div}(\alpha_i) \; +\; \IDp_i - \frac{\deg \IDp_i}{\deg \IDb} \IDb
\; = \; \sum_{j=1}^{\nu} \, c_{ij} \IDa_j \;\; .
\]
The calculation of the $\alpha_i, c_{ij}$ is described in
\cite{PS}.
 
Consequently, the coefficient vectors $(a_1,\cdots,a_{\nu},\lambda)$ 
can be chosen
as $\Z_2$-linear combinations of the rows of the following matrix
$A \in \Z_2^{(\nu+e) \times (\nu +1)}$:\bigskip

\[
 A \; =\; \left( \begin{array}{cccccccc}
2^{n_1} & 0 & \cdots & 0 & 0 & | & 0 \\
0 & 2^{n_2} & \cdots & 0 & 0 & | & 0 \\
\cdot & \cdot & \cdots & \cdot & \cdot  & | & \cdot\\
 \cdot & \cdot & \cdots & \cdot & \cdot & | & \cdot \\
0 & 0 & \cdots & 2^{n_{\nu-1}} & 0 & | & 0 \\
0 & 0 & \cdots & 0 & 2^{n_{\nu}} & | & 0 \\
-- & -- & --- & -- & -- &  & ---\\
   &    &     &    &   & | & \frac{\deg (\IDp_1)}{\deg (\IDb)} \\
   &    & c_{ij} &  &  & | & \vdots \\
   &    &        &  &  & | & \frac{\deg (\IDp_{e})}{\deg (\IDb)} 
\end{array} \right)
\]\bigskip

Each row $(a_1,\cdots,a_{\nu},\lambda)$ of $A$ corresponds to a linear
combination 
satisfying
\begin{equation}
\sum_{j=1}^{\nu} \, a_j \IDa_j + \lambda \IDb \; \equiv \; 
\widetilde{\div}(\alpha ) \, \bmod \; \Dl_F(P\!E) \;\; .
\end{equation}
Condition (\ref{cond2}) gives
\begin{equation}
\prod_{i=2}^m \IDg_i^{b_i} \; =\; \sg (\alpha) \times
\prod_{j=1}^{\nu} \IDe_j ^{a_j} \times
 \IDe_{\IDb}^{\lambda} \;\; .
\end{equation}
Obviously, the family $(\IDg_i) _{2 \leq i \leq m}$ is free over
$\Fq$ implying that the exponents $b_i$ are uniquely defined.
Consequently, if the $k$-th coordinate of the product
$\sg (\alpha) \times \prod_{j=1}^{\nu} \IDe_j ^{a_j}
\times \IDe_{\IDb}^{\lambda}$ is $-1$ we must
have $b_k=1$, otherwise $b_k=0$ for $2 \leq k \leq m$.
(We note that the product over all coordinates is always 1.) 
Therefore, we denote by $b_{2,j},\cdots,
b_{m,j}$ the exponents of the relation belonging to the j-th
column of $A$ for $j=1,\cdots,\nu+e$.

Unfortunately, the elements $\alpha $ are only given up to exceptional units. Hence, we must additionally consider the signs of the  exceptional units of $F$. For
\begin{equation}
\wi{\cal E}{}^{e\!x\!c}_F \; =\; \langle -1 \rangle \times \langle \wi\varepsilon_1,\cdots, \wi\varepsilon_{r+c-1+e} \rangle
\label{utilde1}
\end{equation}
we put:
\begin{equation}
\sg (\widetilde \varepsilon _j ) = \prod_{i=1}^m \IDg_i^{b_{i,j+v+e}} \;\; .
\label{utilde2}
\end{equation}

Using the notations of (\ref{utilde1}) and (\ref{utilde2})  
the rows of the following matrix $A' \in
\Z_2^{(\nu+e+r+c)\times (\nu+m)}$ generate all relations for
the $(\IDa_j,\IDe_j),\; (\IDb,\IDe_{\IDb}),\; ({\bf 0},\IDg_i)$.\bigskip

{\small
\[
 A' \; =\; \left( \begin{array}{ccccccc}
   &  &    & |  & b_{2,1} & \cdots & b_{m,1} \\
   &  &    & |  &  \cdot &  \cdots &  \cdot \\
   & A &   & |  &  \cdot & \cdots & \cdot \\
   &   &   & |  &  \cdot & \cdots & \cdot \\
   &   &   & |  & b_{2, \nu+e} & \cdots & b_{m,\nu+e} \\
-  & --- & - & | &  - &--- & - \\
   &  &  & | & b_{2,\nu+e+1} & \cdots & b_{m,\nu+e+1}  \\
   &  &  & | & \cdot & \cdots & \cdot \\
   &  {\bf O}  &   & | &\cdot & \cdots & \cdot \\
   &   &   & | & \cdot & \cdots & \cdot \\
   &     &   & | & b_{2,\nu+e+r+c} & \cdots & b_{m,\nu+e+r+c}
\end{array} \right) \;\; .
\]
}

\subsection{The algorithm for computing $\wi\Cl{}_F^{\,pos}$ }
 
We assume  that $P\!E=\{\IDp_1,\cdots,\IDp_{e}\} \ne \emptyset$ is ordered by increasing 2-valuations $v_2(\deg \IDp_i)$;  that the group $\Cl_F^{\,pos}$ of positive divisor classes is isomorphic to the direct sum
\[
\Cl_F^{\,pos} \; \cong \;  \oplus_{i=1}^{w} \, \Z / 2^{m_i}\Z \,;
\]
and that we know a full set of representatives $(\IDb_i,\IDf_i)
\;\; (1 \leq i \leq w)$ for all classes. \smallskip

Then each $(\IDb,\IDf)
\in \wi\Dl{}_F^{\,pos}$ satisfies $\deg (\IDb) \in \deg (\Dl_F(P\!E))$
and
\[
\IDb \; \equiv \; \sum_{i=1}^w \, b_i \IDb_i \; 
\bmod (\Dl_F(P\!E)+\widetilde{\Pl}_F) \;\;.
\]
Obviously, we obtain
\[
0 \; \equiv \; \deg (\IDb) \; \equiv \; \sum_{i=1}^w \,
b_i \deg(\IDb_i) \; \bmod \, \deg (\Dl_F(P\!E)) \;\; .
\]
We reorder the $\IDb_i$ if necessary so that
\[
v_2 (\deg (\IDb_1)) \; \leq \; v_2 (\deg (\IDb_i)) \;\;
(2 \leq i \leq w)
\]
is fulfilled. We put
\[
\begin{aligned}
t \; :&= \; \max (\min(\{v_2(\deg (\IDp))\mid \IDp \in \Dl_F(P\!E)\})-v_2(\deg (\IDb_1)),0) \\
&= \; \max (v_2(\deg (\IDp_1))-v_2(\deg (\IDb_1),0)
\end{aligned}
\]
and 
\[
\delta \; :=\; b_1 + \sum_{i=2}^w \, \frac{\deg (\IDb_i)}{\deg (\IDb_1)} b_i
\;\; .
\] 
Then we get:
\[
\IDb \; \equiv \; \sum_{i=2}^w \, b_i \left( \IDb_i - 
\frac{\deg (\IDb_i)}{\deg (\IDb_1)} \IDb_1 \right)
+ \delta \IDb_1 \; \bmod \, (\Dl_F(P\!E)+\widetilde{\Pl}_F) 
\]
and so
\[
\deg \IDb \equiv 0 \equiv \sum b_i \times 0 + \delta \deg \IDb_1 \bmod \deg \Dl_F(P\!E).
\]
From this it is immediate that a full set of 
representatives of the elements of $\wi\Cl{}_F^{\,pos}$
is given by
\[
\left( \IDb_i - \frac{\deg (\IDb_i)}{\deg (\IDb_1)} \IDb_1,
\IDf_i \times \IDf_1^{-\deg (\IDb_i) / \deg (\IDb_1)} \right)
\; \mbox{ for } \; 2 \leq i \leq w
\]
 and
 \[
 (\IDb'_1:= 2^t\IDb_1 - 2^t \frac{\deg \IDb_1}{\deg \IDp_1}\IDp_1,\IDf_1^{2^t}) \;\; .
\]
Let us denote the class of $(\IDc,\IDf)$ in $\wi\Cl{}_F^{\,pos}$
by $[\IDc,\IDf]$.\smallskip

Now we establish a matrix of relations for the generating
classes. For this we consider relations:
\[
\sum_{i=2}^w \, a_i \left[ \IDb_i - 
\frac{\deg (\IDb_i)}{\deg (\IDb_1)} \IDb_1 , \IDf_i \times
\IDf_1^{-\frac{\deg (\IDb_i)}{\deg (\IDb_1)}} \right]
+a_1 \left[ 2^t \IDb'_1, \IDf_1^{2^t} \right] \; =\; 0 \;\; ,
\]
hence
\[
\sum_{i=2}^w \, a_i [ \IDb_i,\IDf_i ] + \left(
2^t a_1 - \sum_{i=2}^w \, \frac{\deg (\IDb_i)}{\deg (\IDb_1)} a_i
\right) [ \IDb_1, \IDf_1 ] \; =\; 0 \;\; .
\]
A system of generators for all relations can then be computed 
analogously to
the previous section. We calculate a basis of the nullspace
of the matrix $A'' =(a''_{ij}) \in \Z^{w \times 2w}$ with first row
\[
\left(2^t, - \frac{\deg (\IDb_2)}{\deg (\IDb_1)}, \cdots , -
\frac{\deg (\IDb_w)}{\deg (\IDb_1)} , 2^{m_1},0, \cdots,0\right)
\]
and in rows $i=2,\cdots,w$ all entries are zero except for
$a''_{ii}=1$ and $a''_{i,w+i}=2^{m_i}$. We note that we are
only interested in the first $w$ coordinates of the
obtained vectors of that nullspace.

\section{Examples}

The methods described here are implemented in the computer algebra system Magma
\cite{C+}.  Many of the fields used in the examples were results of queries
to the QaoS number field database \cite[section 6]{K3}.
More extensive tables of examples can be found at:
\[\mbox{\tt http://www.math.tu-berlin.de/\~{}pauli/K}\]
In the tables abelian groups are given as a list of the orders of their cyclic factors.
\begin{itemize}
\item[] $[:]$ denotes the index $(K_2(O_F):W\!K_2(F))$ (see \cite[equation (6)]{BG});
\item[] $d_F$ denotes the discriminant for a number field $F$;
\item[] $\Cl_F$ denotes the class group, $P$ the set of dyadic places; 
\item[] $\Cl'_F$ denotes the 2-part of $\Cl/\langle P\rangle$;
\item[] $\widetilde{\Cl}_F$ denotes the logarithmic classgroup;
\item[] $\Cl^{\,pos}_F$ denotes the group of positive divisor classes;
\item[] $\wi\Cl{}^{\,pos}_F$ denotes the group of positive divisor classes of degree 0;
\item[] $rk_2$ denotes the 2-rank of the wild kernel $W\!K_2$.
\end{itemize}
\medskip

K. Belabas and H. Gangl have developed an algorithm for the computation of the tame
kernel $K_2\mathcal{O}_F$ \cite{BG}.
The following table contains the structure of $K_2\mathcal{O}_F$ as computed by
Belabas and Gangl and the $2$-rank of the wild kernel $W\!K_2$
calculated with our methods for some imaginary quadratic fields.  
We also give the structure of the wild kernel if it
can be deduced from the structure of $K_2\mathcal{O}_F$ and of the rank of the
wild kernel computed here or in \cite{PS}.

\subsection{Imaginary Quadratic Fields}

\noindent
\begin{tabular}{|r|ccc|cc|cccc|cc|}
\hline
$d_F$    &  $ \Cl_F$  \!\!      & $K_2\mathcal{O}_F$ 
                              & $\![:]\!$ & $\!\!|P|\!\!\!$ & $\!\!\!|P\!E|\!\!\!$ &  $\Cl'_F$  &  $\widetilde{\Cl}^{\phantom{1}}_F$  & $\Cl^{\,pos}_{F_{\phantom{1}}}$\!\! & $\wi\Cl{}_F^{\,pos}$ & $\!rk_2\!\!$ & \!\!$W\!K_2$ \\
\hline
-184 & [ 4 ]       & [ 2 ]    & 1     & 1     &  1     & [ 2 ] & [ 1 ]   & [ 2 ]  & [ ]      & 1               & [ 2 ] \\
-248 & [ 8 ]       & [ 2 ]    & 1     & 1     &  1     & [ 4 ] & [ 2 ]   & [ 4 ]  & [2,2]& 1               & [ 2 ] \\

-399 & [2,8]       & [2,12]   & 2     & 2     & 2      & [ 2 ] & [ 4 ]   & [ 2 ]  & [ 2 ] & 1              & [ 4 ]\\
-632 & [ 8 ]       & [ 2 ]    & 1     & 1     &  1     & [ 4 ] & [ 2 ]   & [ 4 ]  & [2,2]& 1               & [ 2 ] \\

-759 & [2,12]      & [2,18]   & 6     & 2     & 2      & [ 2 ] & [ 2 ] & [ 2 ] & [ 2 ] & 1 & [ 6 ] \\
-799 & [ 16 ]      & [2,4]    & 2     & 2     & 2      & [ 2 ] & [2,4] & [ 2 ] & [ 2 ] & 2 & [2,2] \\

-959 & [ 36 ]      & [2,4] & 2     &  2    &  2     & [ 4 ] & [4,8]& [ 4 ]  & [ 4 ]   & 1               & [ 4 ] \\
\hline
\end{tabular}

\subsection{Real Quadratic Fields}

\noindent
\begin{tabular}{|r|cc|cc|cccc|c|}
\hline
$d_F$ & $\Cl_F$ &[:]      &$\!|P|$\!& $\!\!|P\!E|\!\!$ & $\Cl'$  &  $\widetilde{\Cl}_F$ & $\Cl^{\,pos}_F$ & $\wi{\Cl}{}_F^{{\,pos}^{\phantom{1}}}$ &  $rk_2$\\
\hline
776 & [ 2 ] &  4     & 1 &    1    &    [ 2 ]   &[ ]   &  [ 2,2 ] &       [ 2 ]  & 2 \\
904 & [ 8 ] &  4     & 1 &    1     &    [ 4 ]   &[ 2 ] & [ 4 ]   &[ 2,2 ]       & 1 \\

\hline
29665  & [ 2,16 ] & 8 & 2 & 2 &[ 2 ] & [ 2 ] &[ 2,2 ] & [ 2,2 ] &2 \\
34689  & [ 32 ]   & 8 & 2 & 2 & [ ] & [ ] & [ 2 ] & [ 2 ] & 1 \\
69064  & [ 4,8 ] & 4 & 1 & 1 & [2,8] & [ 8 ] &[ 2,8 ] & [ 8 ] & 2\\
90321  & [2,2,8] & 24& 2 &    2     &    [2,2]&        [2,4]  &   [2,2,2,2]   &    [2,2,2,2] & 4 \\
104584 & [ 4,8 ] & 4 & 1 &    1     &    [2,8]&        [2,4]  &    [ 2,8 ]    &    [ 2,2,4 ]  & 2\\
248584 & [ 4,8 ] & 4 & 1 &    1     &    [2,8]&        [2,4]  &    [2,2,8]    &    [2,2,2,4]  & 3\\
300040 & [2,2,8] & 4 & 1 &    1     &    [2,8]&        [ 8 ]  &    [ 2,8 ]    &    [ 8 ]      & 2 \\
374105 & [ 32 ]  & 8 & 2 &    2     &    [ ]  &         [ ]   &    [ 2 ]      &    [ 2 ]      & 1 \\
\hline
171865 & [ 2,32 ] & 8 & 2 & 2 & [ 4 ] & [ 4 ] & [ 2,2,4 ] & [ 2,2,4 ] & 3\\
285160 & [ 2,32 ] & 4 & 1 & 1 & [ 32 ] & [ 32 ] & [ 32 ] & [ 32 ] & 1\\
318097 & [ 64 ] & 8 & 2 & 2 & [ ] & [ ] & [ 2,2 ] & [ 2,2 ] & 2\\
469221 & [ 64 ] & 12 & 1 & 1 & [ 64 ] & [ 64 ] & [ 2,64 ] & [ 2,64 ] & 2\\
651784 & [ 2,32 ]& 4 & 1 &    1     &    [2,16]&       [2,8]  &   [2,2,16] &   [2,2,2,8] & 3\\
\hline
\end{tabular}

\newpage
\subsection{Examples of Degree 3}

The studied fields are given by a generating polynomial $f$ and have Galois group of their normal closure
isomorphic to $C_3$ (cyclic) or ${\goth S}_3$ (dihedral); $r$ denotes the number of real places.

\begin{center}
\rotatebox{-90}{

\begin{tabular}{|c|rccc|c|cc|cccc|c|}
\hline
&&&&&&&&&&&&\\
$f$ &      $d_F$ & $\!r\!$ & Gal &  $C_F$  & [:]   & $\!|P|\!$ & $\!|P\!E|\!$ &   $\Cl'_F$ & $ \widetilde{\Cl}_F$ & $\Cl^{\,pos}$ & $\widetilde{\Cl}{}^{{\,pos}^{\phantom{1}}}_F$ & $\!rk_2\!$ \\
&&&&&&&&&&&&\\
\hline
&&&&&&&&&&&&\\
$x^3+x^2-18x+12$  &  3957 & 3 &   ${\goth S}_3$ &[ ]   & 48 &               2 & 2 & [ ] &   [ ] &          [ 2 ]   & [ 2 ]    &1\\
$x^3-21x+28$        &  3969 & 3 &   C$_3$ &[ 3 ]& 32 &               3 & 3    & [ ] &   [ ] &          [2,2]& [2,2] &2\\
$x^3-10x+1$         &  3973 & 3 &   ${\goth S}_3$ &[ ]   & 16 &               2 & 2    & [ ] &   [ ] &          [ 2 ]   & [ 2 ]    &1\\
$x^3+x^2-11x-12$  &  3981 & 3 &   ${\goth S}_3$ &[ 2 ]& 16 &               2 & 2 &    [ ] &   [ ] &          [ 2 ]   & [ 2 ]    &1\\
$x^3-16x+4$         &  3988 & 3 &   ${\goth S}_3$ &[ ]   & 8  &               1 & 1    & [ ] &   [ ] &          [ 2 ]   & [ 2 ]    &1\\
&&&&&&&&&&&&\\
\hline
&&&&&&&&&&&&\\
$x^3 - 40x + 1349$        &   -997523 &1&    ${\goth S}_3$ &[ 16 ]   &   4  &               2 & 2     &[ ]        &  [ ]        &        [ 2 ]    &   [ 2 ]    &   1\\
$x^3 - 25x + 198$         &   -996008 &1&    ${\goth S}_3$ &[2,8] &   4   &              2 & 2     &[ 4 ]     &  [ 4 ]      &       [ 2,4 ] &   [ 2,4 ] &   2\\
$x^3 + x^2 - 47x - 1365$  &   -994476 &1&    ${\goth S}_3$ &[ 16 ]   &   6    &             1 & 1     &[ 16 ]    &  [ 16 ]     &       [ 16 ]   &   [ 16 ]   &   1\\
$x^3 + x^2 + 126x + 234$  &   -992696 &1&    ${\goth S}_3$ &[2,8]  &  4     &            2 & 2     &[ 2 ]     &  [ 2 ]      &       [ 2,2 ] &   [ 2,2 ] &   2\\
$x^3 + x^2 + 39x - 155$   &   -992620 &1&    ${\goth S}_3$ &[2,8]   & 2      &           1 & 1     &[ 2,8 ]  &  [ 2,8 ]   &     [ 2,8 ]   & [ 2,8 ]   & 2\\
$x^3 + x^2 + 59x - 63$    &   -991852 &1&    ${\goth S}_3$ &[ 16 ]     & 2       &          1 & 1     &[ 16 ]    &  [ 16 ]     &       [ 16 ]   &   [ 16 ]    &  1\\
$x^3 + x^2 - 108x + 2304$ &   -991423 &1&    ${\goth S}_3$ &[ 16 ]      &8        &         3 & 3     &[ ]        &  [ ]         &       [ 2 ]    &   [ 2 ]     &  1\\
&&&&&&&&&&&&\\
\hline
&&&&&&&&&&&&\\
$x^3 + x^2 - 49x - 48$   & 453317  &3   & ${\goth S}_3$ &[ 16 ]  & 16   &             2  &2     &[ ]  &   [ ]   &      [ 2 ]     &     [ 2 ]    &      1\\
$x^3 - 203x + 548$      &  1014140 &3  &  ${\goth S}_3$ &[ 16 ] &  32  &              2 & 1     &[ ] &    [ ]  &       [ 2 ]    &      [ 2 ]   &       1\\
$x^3 + x^2 - 164x + 64$   &1085681 &3   & ${\goth S}_3$ &[ 16 ]  & 32    &            3  &3     &[ 2 ] & [ 2 ] &[2,2,2,2]&[2,2,2,2]&4\\
&&&&&&&&&&&&\\
\hline
&&&&&&&&&&&&\\
$x^3 - x + 216$           & -314927 &1 &   ${\goth S}_3$ &[ 64 ]    &8        &         3  &3     &[ ]    &[ ]      &[ 2 ]    &[ 2 ]    &1\\
$x^3 + x^2 - 232x - 1840$ &-526836 &1 &   ${\goth S}_3$ &[2,32] &12    &            2  &2     &[ 2 ] &[ 2 ] &  [ 2,2 ] &[ 2,2 ] &2\\
$x^3 + 70x + 236$      &   -718948 &1 &   ${\goth S}_3$ &[ 64 ]   & 8    &             2  &1     &[ 2 ] &[ 2 ] &  [ 2 ]    &[ 2 ]    &1\\
&&&&&&&&&&&&\\
\hline
\end{tabular}
}
\end{center}

\subsection{Examples of Higher Degree}

\begin{center}
\rotatebox{-90}{
\begin{tabular}{|c|rccc|c|cc|cccc|c|}
\hline
&&&&&&&&&&&&\\
$f$ &      $d_F$ & $\!r\!$ & Gal  & $\Cl_F$  & [:]   & $\!|P|\!$ & $\!|P\!E|\!$ &  $\Cl'_F$ & $ \widetilde{\Cl}_F$ & $\Cl^{\,pos}_F$ & $\wi{\Cl}{}_F^{{\,pos}^{\phantom{1}}}$ & $\!rk_2\!$ \\
&&&&&&&&&&&&\\
\hline
&&&&&&&&&&&&\\
$x^4 - 59x^2 - 120x - 416$   &    -860400 &2 &   D$_4$ &[ 16 ] &     8  &         2 & 2 &    [ ] &      [ ]  &        [ 2 ]  &  [ 2 ]   & 1\\
$x^4 - x^3 - 2x^2 + 5x + 1$   &   -3967 &2   & ${\goth S}_4$ &[ ] &   8        &         2  &2   &   [ ]  &  [ ]  &   [ 2 ]   &  [ 2 ]   &1 \\
$x^4 - x^3 + 86x^2 - 66x + 1791$ &701125 & 0  &  D$_4$ &[ 2,8 ] &   1    &       1  &1     &[ 2,8 ] &[ 2,8 ] &[ 2,8 ] &[ 2,8 ] &2\\
$x^4 + 14$        &                 702464 & 0   & D$_4$ &[ 4,4 ] &   1     &      1  &1 &    [ 4 ] &   [ 2 ]  &[ 4 ]    &[ 2,2 ] &1\\
$x^4 + 58x^2 + 1$ &                  705600 &0   & E$_4$& [ 4,8 ] & 2       &     2  &2  &   [ 4 ]  &  [ 2 ]  & [ 4 ]  &  [ 2,2 ] &1\\
$x^4 - 2x^3 + 59x^2 - 24x + 738$  &728128 &0 &   D$_4$ &[ 32 ]  &  2       &    2  &2   &  [ 2 ]   & [ 4 ]  &[ 2 ]  &  [ 2 ]    &1\\
$x^4 + 21x^2 + 120$   &              730080 &0  &  D$_4$ &[ 4,8 ] & 6         &   2  &2    & [ 2 ]    &[ 4 ]  & [ 2,2 ] &[ 2,2 ] &2\\
$x^4 - 5x + 30$        &             766125 &0   & ${\goth S}_4$ &[ 2,16 ] &20         &  3  &3 &    [ ]       &[ ]     & [ ]     &  [ ]      & 0\\
$x^4 + 58x^2 + 1$       &            705600 &0    &E$_4$ &[ 4,8 ] & 2           & 2  &2  &   [ 4 ]    &[ 2 ] &  [ 4 ] &   [ 2,2 ] &1\\
&&&&&&&&&&&&\\
\hline
&&&&&&&&&&&&\\
$x^5 + x^4 + x^3 - 8x^2 - 12x + 16$      & -4424116 &3   & ${\goth S}_5$ &[ 4 ]  &64 &   3  &2   &  [ ] & [ ] &   [ 2 ]  &[ 2 ]   &1\\
$x^5 + x^4 - 13x^3 - 26x^2 - 8x - 1$  &-3504168 &3  &  ${\goth S}_5$ &[ 4 ] &16    &            2 & 2 &   [ ]   & [ ]   &  [ 2,2 ] &[ 2,2 ] &2\\
$x^5 - 10x^3 + 9x^2 + 7x - 1$        & -3477048 &3  &  ${\goth S}_5$ &[ 4 ] &16  &              2 & 2 &   [ ]   & [ ]   &     [ ]   &    [ ]   &    0\\
$x^5 + 2x^4 + 6x^3 + 11x^2 - 2x - 9$ &-3420711 &3&    ${\goth S}_5$ &[ 4 ] &8&                 1  &1  &  [ 4 ] &[ 4 ] & [ 4 ]   & [ 4 ]    &1\\
$x^5 - 14x^3 + 26x^2 - 11x - 1$       &-3356683 &3  &  ${\goth S}_5$ &[ 4 ] &16 &              2  &2  &  [ ]    &[ ]    &    [ 2 ]  &  [ 2 ]  &  1\\
$x^5 + 2x^4 + 9x^3 + 3x^2 + 10x - 24$  & 2761273  &1    &${\goth S}_5$ &[ 10 ] &8   &  3  &3    & [ ] & [ ]   & [ ]    & [ ]  &    0\\
$x^5 + x^4 - 3x^3 + 15x^2 + 36x - 18$   & 3825936  &1 &   D$_5$ &[ 11 ] &288  & 3  &1     &[ ] & [ ]  & [ ]    & [ ]   &   0\\
$x^5 + 2x^4 + 12x^3 + 14x^2 - 12x - 16$ &4892116  &1  &  ${\goth S}_5$ &[ 12 ] &8 &    3  &3  &   [ ] & [ ]  &    [ 2 ] & [ 2 ]  & 1\\
$x^5 + 2x^4 - 8x^3 - 4x^2 + 7x + 1$  &   13664837 &5 &   ${\goth S}_5$ &[ 4 ] & 64 &   2  &2 &    [ ] & [ ]   &    [ 2 ] & [ 2 ] &  1\\
$x^5 + 2x^4 - 11x^3 - 27x^2 - 10x + 1$ & 17371748 &5  &  ${\goth S}_5$ &[ 2 ]  &64    &2  &2     &[ ] & [ ]  & [ 2 ] & [ 2 ]  & 1\\
&&&&&&&&&&&&\\
\hline

\end{tabular}
}
\end{center}


\parbox[t]{6,2cm}{\address{Jean-Fran\c{c}ois {\sc Jaulent}\\
Universit\'e Bordeaux I\\
Institut de Math\'ematiques\\
351, Cours de la Lib\'eration\\
33405 Talence Cedex, France\\
{\small \tt jaulent@math.u-bordeaux1.fr}
}}
\parbox[t]{6,4cm}{\address{Sebastian {\sc Pauli}\\
University of North Carolina\\
Department of Mathematics and \\Statistics\\ 
Greensboro, NC 27402, USA\\
{\small \tt s\_pauli@uncg.edu}}
}\medskip

\parbox[t]{6,2cm}{\address{Michael E. {\sc Pohst}\\
Technische Universit\"at Berlin\\
Institut f\"ur Mathematik   MA 8-1\\
Stra{\ss}e des 17. Juni 136\\
10623 Berlin, Germany\\
{\small \tt pohst@math.tu-berlin.de}
}}
\parbox[t]{6,4cm}{\address{Florence {\sc Soriano-Gafiuk}\\
Universit\'e de Metz\\
D\'epartement de Math\'ematiques\\
Ile du Saulcy\\
57000 Metz, France\\
{\small \tt soriano@poncelet.univ-metz.fr}
}}
\end{document}